\magnification=1100
\input amssym.def
\input amssym.tex
 \font\newrm =cmr10 at 24pt
\def\bul{\raise .9pt\hbox{\newrm .\kern-.105em } }

 \def\fr{\frak}

 \def\h{\hbox{ }}

 \def\u{{\fr u}}
 \def\r{{\fr r}}
 
 \def\m{{\fr m}}
 \def\n{{\fr n}}
 \def\a{{\fr a}}

 \def\ss{{\fr s}}
 \def\k{{\fr k}}
 \def\b{{\fr b}}
 
 \def\hh{{\fr h}}
 \def\tt{{\fr t}}

 \def\g{{\fr g}}

 \def\q{{\fr q}}

 \def\<{\le}
 \def\>{\ge}

 \def\s{{\h\subset\h}}
 
 \def\vs{\vskip }

 \def\mapright#1
  {\smash{\mathop
  {\longrightarrow}
  \limits^{#1}}}

 \def\kk#1{{\kern .4 em} #1}
 \def\vs{\vskip 1pc}

\font\twelverm=cmr12 at 14pt
\font\ninerm=cmr10 at 9pt
\hsize = 31pc
\vsize = 45pc
\overfullrule = 0pt
\font\smallbf cmr10 at 9pt
\font\authorfont=cmr10 at 10pt
\font\ninerm=cmr10 at 9pt
\font\smallrm=cmr10 at 7pt
\centerline{\twelverm Root Systems for Levi Factors and }
\vskip 4pt
\centerline {\twelverm Borel--de Siebenthal Theory} 
\vskip 1.5pc
\baselineskip=11pt
\vskip8pt
\centerline{\authorfont BERTRAM KOSTANT\footnote*{\ninerm Research supported in part by the KG\&G
Foundation}}\vskip 2pc
\baselineskip 20pt
\hfill \it Dedicated to Gerry Schwarz on the occasion of his 60th birthday
\vskip 2pc
\baselineskip=14pt 
\noindent{\smallbf ABSTRACT. \ninerm $\,$ Let $\m$ be a Levi factor of a
proper parabolic subalgebra $\q$ of a complex semi\-simple Lie algebra $\g$. Let $\tt = \hbox{\rm cent}\,\m$.
A  nonzero element
$\nu\in \tt^*$ is called a $\tt$-root if the corresponding adjoint weight space $\g_{\nu}$ is not zero. If
$\nu$ is a $\tt$-root, some time ago we  proved that $\g_{\nu}$ is $\hbox{\rm ad}\,\m$ irreducible. Based on
this result  we develop in the present paper a theory of $\tt$-roots which replicates much of the structure of
classical root theory (case where $\tt$ is a Cartan subalgebra). The results are applied to obtain new results
about the structure of the nilradical $\n$ of $\q$. Also applications in the case where $\hbox{\rm dim}\,\tt=1$
are used in Borel--de Siebenthal theory to determine irreducibility theorems for certain equal rank subalgebras
of $\g$. In fact the irreducibility results readily yield a proof of the main assertions of the Borel--de
Siebenthal theory.}

\vskip 1pc
\baselineskip 16pt

\centerline{\bf 0. Introduction}
 \rm \vskip 1.5pc {\bf 0.1.} \rm Let $\m$ be a Levi factor of a
proper parabolic subalgebra $\q$ of a complex semisimple Lie algebra $\g$. Let $\tt = \hbox{\rm cent}\,\m$
and let $\ss = [\m,\m]$ so that one has a direct sum $\m = \tt + \ss$. Let $\r$ be the Killing form
orthocomplement of $\m$ in $\g$ so that $\g = \m + \r$ and $[\m,\r]\s \r$. A nonzero element
$\nu\in \tt^*$ is  called a $\tt$-root if $\g_{\nu}\neq 0$ where $\g_{\nu} = \{ z\in \g\mid \hbox{\rm ad}\,x (z)
=
\nu(x)\,z,\,\,\forall x\in \tt\}$. One readily has $\g_{\nu}\s \r$ and a direct sum $$\r =
\sum_{\nu\in R}\,\,\g_{\nu}\eqno (0.1)$$ where $R\s \tt^*$ is the set of all $\tt$-roots. It is immediate
that $\g_{\nu}$ is an $\hbox{\rm ad}\,\m$-submodule of $\r$ for any $\nu\in R$. Some time ago we proved \vs {\bf
Theorem 0.1.} {\it $\g_{\nu}$ is an irreducible $\hbox{\rm ad}\,\m$-module for any $\nu\in R$ and any irreducible
$\m$-submodule of $\r\,$ is of this form. In particular $\r$ is a multiplicity-free $\hbox{\rm ad}\,\m$-module
and (0.1) is the unique decomposition of $\r$ as a sum of irreducible $\hbox{\rm ad}\,\m$-modules.}

\vs Our
Theorem 0.1 appeared, with the apppropriate citations, as Theorem 8.13.3 in [Wol] and Theorem 2.1 in [Jos]. In
the present paper we will use Theorem 0.1 (reproved for convenience) to develop a theory of $\tt$-roots which in
many ways replicates results in the usual root theory, i.e., the case where $\tt$ is a Cartan subalgebra of
$\g$. For example it is established that if $\mu,\nu\in R$ and $\mu + \nu\in R$,  then one has the
equality $$[\g_{\mu},\g_{\nu}] = \g_{\mu +\nu}.$$ Also with respect to a natural inner product on $\tt^*$
if $\mu,\nu\in R$ and $(\mu,\nu)<0$, then $\mu + \nu\in R$,  and if $(\mu,\nu) >0$,  then $\mu-\nu\in R$. (See
Theorem 2.2.)

The nilradical $\n$ of $\q$ is contained in $\r$, and one introduces a set $R^+$ of positive $\tt$-roots so
that $$\n= \sum_{\nu\in R^+}\,\,\g_{\nu}.$$ As in the Cartan subalgebra case one can similarly define the
set $R_{\hbox{\smallrm simp}}\s R^+$ of simple positive $\tt$-roots and prove that if, by definition, $\ell(\tt)=
\hbox{\rm dim}\,
\tt$, then
$\hbox{\rm card}\, R_{\hbox{\smallrm simp}} = \ell(\tt)$. (See Theorem 2.7.) In fact if $R_{\hbox{\smallrm simp}}
=
\{\beta_1,\ldots,\beta_{\ell(\tt)}\}$,   then the
$\beta_i$ are a basis of $\tt^*$ and $(\beta_i,\beta_j)\leq 0$ if $i\neq j$. In addition one proves
  that $\n$ is
generated by $\g_{\beta_i}$ for $i=1,\ldots,\ell(\tt)$. In fact we prove that for the nilradical $\n$ of a
parabolic subalgebra of $\g$ one has

 \vs {\bf Theorem 0.2.} {\it Except for indexing, the upper central
series of $\n$ is the same as the lower central series of $\n$.}\vs In Section 3 of the paper we
deal with the case where $\ell(\tt) =1$ so that $\q$ is a maximal parabolic subalgebra. The results are applied
here to Borel--de Siebenthal theory and irreducibility theorems are obtained for the adjoint action of equal (to
that of $\g$) rank subalgebras $\g^{a_j}$ of $\g$ on the Killing form orthocomplement of $\g^{a_j}$ in $\g$. 
 In Remark 3.9 we also show that these results provide a proof of the main statements of the Borel--de Siebenthal
theory. 

\vskip 1.5 pc
\centerline{\bf 1. Levi factor root system}\vskip 1pc
{\bf 1.1.} Let $\g$ be a complex semisimple Lie algebra. Let $\hh$
be a Cartan subalgebra of $\g$ and let $\ell$ be the rank of $\g$.
Let
$\hh^*$ be the dual space to $\hh$ and let $\Delta\s \hh^*$ be the
set of roots for the pair $(\hh,\g)$. For each $\varphi\in \Delta$,
let
$e_{\varphi}\in \g$ be a corresponding root vector. Let $\Delta_+\s
\Delta$ be a choice of positive roots and let $\Pi\s \Delta$ be the
set of simple positive roots. Let $\Delta_-=-\Delta_+$. For any
$\varphi\in
\Delta
$ and
$\alpha\in \Pi$, let $n_{\alpha}(\varphi)\in \Bbb Z$ be the integer
(nonnegative if $\varphi\in \Delta_+$ and nonpositive if $\varphi\in
\Delta_-)$ so that $$\varphi =\sum_{\alpha\in
\Pi}\,n_{\alpha}(\varphi)\,\alpha. \eqno (1.1)$$  If $\u\s \g$ is
any subspace which is stable under $\hbox{\rm ad}\,\hh$, let
$$\eqalign{\Delta(\u)&=\{\varphi\in \Delta\mid e_{\varphi}\in \u\}
\cr \Delta_+(\u)&= \Delta(\u)\cap\Delta_+\cr
\Delta_-(\u)&= \Delta(\u)\cap\Delta_-.\cr}\eqno (1.2)$$\vskip .5pc
Let
$\b$ be the Borel subalgebra of $\g$, containing $\hh$ such that $\Delta(\b)
= \Delta_+$. Let $\n_{\b} = [\b,\b]$ be the nilradical of $\b$. A
standard parabolic subalgebra
$\q$ of
$\g$ is any Lie subalgebra of $\g$ which contains $\b$. 

Let $B$ be the Killing form on $\g$. Assume that $\q$ is some fixed 
standard parabolic
subalgebra of
$\g$. Then $\q$ admits a unique Levi decomposition $$\q = \m +
\n\eqno (1.3)$$ where $\n= \n_{\q}\,(\subset \n_{\b})$ is the nilradical
of
$\q$ and $\m = \m_{\q}$ is the unique Levi factor of $\q$ which
contains $\hh$. We will assume throughout that $\q \neq \g$ so that $\n \neq 0$. Let $\ss =
[\m,\m]$ so that
$\ss$ is the unique maximal semisimple ideal in $\m$. Let $\tt$ be the center of $\m$ so
that $B$ is nonsingular on both $\tt$ and $\m$ and $$\m = \tt
\oplus \ss\eqno (1.4)$$ is a $B$-orthogonal decomposition of $\m$ into a
direct sum of ideals. Let $\hh(\ss) = \hh\cap \ss$ so that $\hh(\ss)$ is
a Cartan subalgebra of $\ss$.
$B$ is nonsingular on $\hh(\ss)$ and $$\hh = \tt \oplus \hh(\ss)
\eqno (1.5)$$
is a $B$-orthogonal decomposition. 

The nonsingular bilinear form, $B|\hh$, on $\hh$ induces a nonsingular
bilinear form on $\hh^*$ which we denote by $B|\hh^*$. We may embed the
dual spaces $\tt^*$ and $\hh(\ss)^*$ to $\tt$ and $\hh(\ss)$,
respectively, in $\hh^*$ so that $\tt^*$ is the orthocomplement of
$\hh(\ss)$ and
$\hh(\ss)^*$ is the orthocomplement to $\tt$. Then $B|\hh^*$ is
nonsingular on both $\tt^*$ and $\hh(\ss)^*$ and $$\hh^* = \tt^* \oplus
\hh(\ss)^*\eqno (1.6)$$ is a $B|\hh^*$ orthogonal direct sum. 
 
Let $\hh_{\Bbb R}^*$ be the real form of $\hh^*$ spanned over $\Bbb
R$ by
$\Delta$. As one knows, $B|\hh^*$ is positive definite on $\hh_{\Bbb
R}^*$. On the other hand, similarly, $\Delta(\ss)$ clearly spans over
$\Bbb R$, a real form, $\hh(\ss)^*_{\Bbb R}$, of $\hh(\ss)^*$. Of course 
$\hh(\ss)^*_{\Bbb R}$ is a real subspace of $\hh_{\Bbb R}^*$ and,
clearly, if $\tt^*_{\Bbb R}$ is the $B|\hh^*$ orthocomplement of
$\hh(\ss)^*_{\Bbb R}$ in $\hh_{\Bbb R}^*$, then $\tt^*_{\Bbb R}$ is a
real form of $\tt^*$ and $$\hh_{\Bbb R}^* = \tt^*_{\Bbb R}\oplus
\hh(\ss)^*_{\Bbb R}\eqno (1.7)$$ is a real Hilbert space orthogonal
direct sum. For any $\gamma\in \hh_{\Bbb R}^* $ we let $\gamma_{\tt}\in
\tt^*_{\Bbb R}$ and
$\gamma_{\ss}\in \hh(\ss)^*_{\Bbb R}$, respectively, be the components of
$\gamma$ with respect to the decomposition (1.7) so that $$\gamma = \gamma_{\tt}
+\gamma_{\ss}\eqno (1.8)$$ and $$(\gamma_{\tt}, \gamma_{\ss})= 0\eqno (1.9)$$
where $(\mu,\lambda)$
denotes the $B|\hh^*$-pairing of any $\mu,\lambda\in \hh^* $. 

Let ${\overline \n}$ be the span of all $e_{-\varphi}$, for $\varphi\in
\Delta(\n)$ so that one has a triangular decomposition $$\g= \m +\n
+{\overline \n}.\eqno (1.10)$$ Now put $\r = \n
+{\overline \n}$ so that $\r$ is $\hbox{\rm ad}\,\m$-stable and one has a
$B$-orthogonal decomposition $$\g = \m + \r.\eqno (1.11)$$\vskip 4pt {\bf
Remark 1.1.} From the general properties of Levi factors  of
parabolic subalgebras one knows that $\m$ is the centralizer of $\tt$
in
$\g$ so that
$$\Delta(\r) =
\{\varphi\in \Delta\mid \varphi_{\tt} \neq 0\}.\eqno (1.12)$$ 

\vskip
.5pc{\bf 1.2}. Let $V$ be a $\tt$-module
and let
$\mu\in
\tt^*$. Put $$V_{\mu}=\{ v\in V\mid x\cdot v = \langle
\mu,x \rangle\,\,v,\,\,\forall x\in \tt\}. \eqno (1.13)$$ The
subspace
$V_{\mu}$ is called the $\mu$-weight space (for $\tt$) of $V$. If
$V_{\mu}\neq 0$, then $\mu$ is called a $\tt$-weight of $V$ and
any
$v\in V_{\mu}$ is called a $\mu$-weight vector. 

If
$V$ is a finite-dimensional $\g$-module, then necessarily $\gamma\in
\hh^*_{\Bbb R}$ where $\gamma$ is any
$\hh$-weight of $V$. One then notes $\mu\in \tt^*$ is a
$\tt$-weight of $V$ if and only if $$\mu = \gamma_{\tt}\,\,\hbox{where}\,\,\gamma\,\,\hbox{is an $\hh$-weight of
$V$}.\eqno(1.14)$$
	An important
special case is when $V = \g$ and the module structure is defined by
the adjoint action. Let $R'$ be the set of all $\tt$-weights of $\g$.
If $V$ is any $\g$-module and $\xi$ is a $\tt$-weight of $V$, it is
obvious that, for any $\mu\in R'$, $$\g_{\mu}\cdot V_{\xi}\s V_{\mu
+\xi.}\eqno (1.15)$$ Clearly $0\in R'$ and $$\g_0 = \m\eqno (1.16)$$
so that $V_{\xi}$ is an $\m$-module.

\vs {\bf Remark 1.2.} If $V$ is
finite dimensional then, since $\ss\s \m$ and $\ss$ is semisimple, note
that
$V_{\xi}$ is a completely reducible $\ss$-module.\vs If $V$ is
equal to the adjoint $\g$-module $\g$ in (1.15) one has
$$[\g_{\nu},\g_{\mu}]\s \g_{\nu + \mu}\eqno (1.17)$$ for any
$\nu,\mu\in R'$. In particular, if $\mu\in R'$, then $$[\m,\g_{\mu}]\s
\g_{\mu}\eqno (1.18)$$ and $$\g_{\mu}\,\,\hbox{is a completely
reducible $\ss$-module for the maximal semisimple ideal $\ss$ of
$\m$.}\eqno (1.19)$$

Let $R = R'\setminus\{0\}$ so that, recalling Remark 1.1, $$R= \{\nu\in
\tt^*\mid \nu = \varphi_{\tt}, \,\,\hbox{for some}\,\, \varphi\in
\Delta(\r)\},\eqno (1.20)$$ and one readily has (see (1.16)) the
direct sum $$\g = \m +
\sum_{\nu\in R}\, \g_{\nu}.\eqno (1.21)$$ We refer to the elements in
$R$ as $\tt$-roots (in $\g$) and $\g_{\nu}$ as the $\tt$-root space 
corresponding to $\nu\in R$. Partially summarizing one readily has

\vs {\bf Proposition 1.3.} {\it One has the two disjoint unions
$$\eqalign{\Delta &= \sqcup_{\mu\in R'}\, \Delta(\g_{\mu})\cr
\Delta(\r)&=\sqcup_{\nu \in R}\,\Delta(\g_{\nu}).\cr}\eqno (1.22)$$
Furthermore if $\nu\in R$, then $$\eqalign{\Delta(\g_{\nu})&=\{\varphi\in
\Delta\mid \varphi_{\tt} = \nu\}\,\,\hbox{and}\cr \{e_{\varphi}\mid
\varphi&\in
\Delta(\g_{\nu}),\,\hbox{is a basis of
$\g_{\nu}$}\}.\cr}\eqno(1.23)$$}\vskip .5pc Let 
$(x,y)$ denote the pairing of
$x,y\in
\g$ defined by
$B$.

\vs {\bf Remark 1.4.} Note that
$\nu\in R$ if and only if $-\nu\in R$ and $$\Delta(\g_{-\nu}) = -
\Delta(\g_{\nu}).\eqno (1.24)$$ Furthermore if $\mu,\nu\in R$ and
$\nu\neq -\mu$, then
$$(\g_{\mu},\g_{\nu})= 0\eqno (1.25)$$ and
$$\g_{\nu}\,\,\hbox{and}\,\,\g_{-\nu}\,\,\hbox{are nonsingularly paired
by
$B$}.\eqno (1.26)$$\vskip .5pc {\bf 1.3.} Let $\tau:\hh\to \hh^*$ be
the linear isomorphism defined by $B|\hh$. Thus for $x\in \hh$ and
$\mu\in \hh^*$, $$\eqalign{\langle \mu,x\rangle&=
(\mu,\tau (x))\cr&=(\tau^{-1}(\mu),x).\cr}\eqno (1.27)$$ Thus if
$\tt_{\Bbb R} = \tau^{-1}(\tt^*_{\Bbb R})$, then one readily has 

\vs {\bf Proposition 1.5.} {\it $$\eqalign{(1)&\,\,\tt_{\Bbb
R}\,\,\hbox{is a real form of}\,\, \tt\cr(2)&\,\,B\,\,\hbox{is real
and positive definite on
$\hh_{\Bbb R}$}\cr(3)&\,\,\tt_{\Bbb R}=\{x\in \tt\mid \langle
\nu,x\rangle\,\,\hbox{is real}\,\,\forall \nu\in R\}.\cr}\eqno (1.28)
$$}\vskip .5pc Let $\nu\in R$. Clearly $\hbox{\rm Ker}\,\nu$ has codimension 1 in
$\tt$. It follows from Proposition 1.5 that there exists a unique
element $h_{\nu}\in \tt_{\Bbb R}$ which is $B$-orthogonal to $\hbox{\rm Ker}\,\nu$
and such that $$\langle \nu,h_{\nu}\rangle = 2.\eqno (1.29)$$ Note
that it follows from (1.27) that $$\tau(h_{\tau}) =
2\,\nu/(\nu,\nu).\eqno (1.30)$$

\vskip .5pc Let $\nu\in R$ and put
$\m(\nu) = [\g_{\nu},\g_{-\nu}]$ so that $\m(\nu)$ is an ideal of
$\m$. Let $\tt(\nu) = \m(\nu) \cap \tt$ and $\ss(\nu) = \m(\nu) \cap
\ss$. Decomposing the adjoint action representation of $\ss$ on
$\m(\nu)$ into its primary components, it is clear that $\tt(\nu)$ is
the primary component corresponding to the trivial representation and
$\ss(\nu)$ is the sum of the remaining components. Consequently one
has the direct sum $$\m(\nu) =  \tt(\nu) \oplus \ss(\nu).\eqno (1.31)$$
\vskip .5pc 

{\bf Remark 1.6.} If $\ss_i$ is a simple component of
$\ss$, one has $\Delta(\ss_i) = -\Delta(\ss_i)$ so that $B|\ss_i$ is
nonzero and hence nonsingular by simplicity. On the other hand if
$\ss_i$ and
$\ss_j$ are distinct simple components then, of course,
$[\ss_i,\ss_j] = 0$. Consequently $\ss_i$ and $\ss_j$ are 
$B$-orthogonal by the invariance of $B$ and the equality
$[\ss_i,\ss_i] = \ss_i$.

\vs {\bf Proposition 1.7.} {\it  Let
$\nu\in R$. Then $$\tt(\nu) = \Bbb C\,h_{\nu}.\eqno (1.32)$$ 
In addition $B|\m(\nu)$ is nonsingular and the kernel of the
adjoint action of
$\m$ on
$\g_{\nu}$ is the orthocomplement of $\m(\nu)$ in $\m$. In particular
$\m(\nu)$ operates faithfully on $\g_{\nu}$.}\vs 

{\bf Proof.} Let
$\m(\nu)^{\perp}$ be the $B$-orthogonal subspace to $\m(\nu)$
in $\m$. Let
$x\in \m, \,y\in \g_{\nu}$,  and $z\in \g_{-\nu}$. But $(x,[y,z]) =
([x,y],z)$. Since $\g_{\nu}$ and $\g_{-\nu}$ are orthogonally paired
by $B$ this proves that $\m(\nu)^{\perp}$ is the kernel of the
adjoint action of $\m$ on $\g_{\nu}$. But then 
$\hbox{\rm Ker}\,\nu= \m(\nu)^{\perp}\cap \tt$. Recalling the definition and properties of $h_{\nu}$ (see (1.29)
and Proposition 1.5) it follows immediately that $B|\hbox{\rm Ker}\,\nu$ is nonsingular. It then follows that
$\tt(\nu)$ must be the one-dimensional
$B$-orthocomplement of 
$\hbox{\rm Ker}\,\nu$ in $\tt$.  But then one has (1.32) by definition of $h_{\nu}$. But now $\ss(\nu)$ is
clearly an ideal in $\ss$. Hence $\ss(\nu)$ is a sum of simple components of
$\ss$. Thus $B|\m(\nu)$ is nonsingular by Remark 1.6. QED

\vs {\bf
1.3.} In this section we will mainly be concerned with decomposing
$\r$ into irreducible $\m$-modules. Effectively this comes down to
understanding the action of $\ss$ on $\g_{\nu}$ for any $\nu\in R$. 
 
Obviously $\Delta_+(\ss)$ defines a choice of positive roots in
$\Delta(\ss)$ so that $$\b(\ss) = \hh(\ss) + \sum_{\varphi\in
\Delta_+(\ss)}\,\,{\Bbb C}\,e_{\varphi}$$ is a Borel subalgebra of
$\ss$. Highest weights and highest weight vectors for $\ss$-modules
will be defined with respect to $\b(\ss)$. 

Let $C(\ss)\s \hh(\ss)^*_{\Bbb R}$ be the dominant Weyl chamber.

\vs
{\bf Proposition 1.8.} {\it Let $\xi,\eta\in C(\ss)$. Then
$$(\xi,\eta) \geq 0.\eqno (1.33)$$}
\indent {\bf Proof.} Let $I$ be an
index set for the simple components of $\ss$ where, if $i\in I$, then
$\ss_i$ is the  corresponding component. One readily has $$C(\ss) =
\sum_{i\in I} C(\ss_i)$$ where $\hh(\ss_i)=
\hh\cap \ss_i$ and $C(\ss_i)\s \tau(\hh(\ss_i))$ is the dominant Weyl
chamber for $\ss_i$. But then, if $i,j\in I$ are
distinct, $C(\ss_i)$ and
$C(\ss_j)$ are $B|\hh^*$-orthogonal by (1.27) and Remark 1.6. Thus it
suffices to prove that, if $i\in I$, and $\xi,\eta\in C(\ss_i)$ are
nonzero, then $$(\xi,\eta)>0. \eqno (1.34)$$ But $B|\ss_i$ is a positive
multiple of the $\ss_i$-Killing form and for the $\ss_i$-Killing form (1.34) is known
(see e.g., Lemma 2.4 in [Kos]). QED

\vs We established the
following theorem some time ago. It appears in the literature in works of
J.~Wolf and A.~Joseph with proper citations in both cases. See Theorem
8.13.3 in [Wol] and in a closer reproduction of my argument, Theorem 2.1 in [Jos].

\vs {\bf Theorem 1.9.} {\it Let $\nu\in R$. Then the
$\tt$-root subspace 
$\g_{\nu}$ is an irreducible
$\m$ and irreducible $\ss$-module under the adjoint action. In fact it
is a faithful irreducible $[\g_{\nu},\g_{-\nu}]\s \m$-module and in
the notation of (1.31) an irreducible $\ss(\nu)$-module. In addition
$$\r = \sum_{\nu\in R} \,\g_{\nu}\eqno (1.35)$$ is a multiplicity-one
representation of $\m$ and the summands (i.e., the $\tt$-root spaces) on
the right hand side of (1.35) are the irreducible components.}\vs 

{\bf
Proof.} Each of the summands on the right side of (1.35) affords a
different character of $\tt$ and hence these summands are inequivalent
as $\m$-modules. Recalling Proposition 1.7 it suffices only to prove
that $\g_{\nu}$ is an irreducible $\ss$-module. The elements
$\{e_{\varphi}\mid \varphi\in \Delta(\g_{\nu})\}$ are a weight basis of
$\g_{\nu}$ for the Cartan subalgebra $\hh(\ss)$. Moreover,  since root
spaces for $\hh$ have multiplicity-one,  the $\hh(\ss)$-weights in
$\g_{\nu}$ have multiplicity-one since $\varphi_{\tt} = \varphi'_{\tt}$
for $\varphi,\varphi'\in \Delta(\g_{\nu})$.

 Assume $\g_{\nu}$ is not
$\ss$-irreducible. Then there exists distinct $\varphi,
\varphi'\in\Delta(\g_{\nu})$ such that $e_{\varphi}$ and $e_{\varphi'}$
are $\ss$-highest weight vectors. In particular $\varphi_{\ss}$ and
$\varphi'_{\ss}$ are in $C(\ss)$. But then
$$(\varphi_{\ss},\varphi'_{\ss})\geq 0\eqno (1.36)$$ by Proposion 1.8.
But $$\eqalign{\varphi &= \nu + \varphi_{\ss}\cr
\varphi' &= \nu + \varphi'_{\ss}.\cr} \eqno (1.37)$$ Hence
$$(\varphi,\varphi')>0.\eqno (1.38)$$ Thus $\beta = \varphi-\varphi'$
is a root. Furthermore $\beta_{\tt}= 0$ so that $\beta\in
\Delta(\m)= \Delta(\ss)$.  Without loss of generality we may choose the ordering so
that $\beta\in \Delta_+(\ss)$. But then $[e_{\beta},e_{\varphi'}]$ is a
nonzero multiple of $e_{\varphi}$. This contradicts the fact that
$e_{\varphi'}$ is an $\ss$-highest weight vector. QED

\vskip 1.5 pc 
\centerline {\bf 2. Properties of the $\tt$-root system}\vskip 1pc
{\bf 2.1.} We will utilize Theorem 1.9 to establish some properties
of $R$. To begin with 

\vs {\bf Lemma 2.1.} {\it Assume $\nu,\mu\in R$
and $\nu +\mu\neq 0$. Assume also that $[\g_{\nu},\g_{\mu}] \neq 0$.
Then $\nu + \mu \in R$ (obvious)  
and one has the equality $$[\g_{\nu},\g_{\mu}] = \g_{\nu + \mu.}\eqno
(2.1)$$}

 {\bf Proof.} The left side of (2.1) is a nonzero
$\m$-submodule of the right side. One therefore has the equality (2.1)
by irreducibility. QED

\vs Let $p,q\in \Bbb Z$ where $p\leq q$. Let
$ I_{p,q}$ denote the set of integers $m$ such that  $p\leq m\leq q$. A finite
nonempty subset $I\s \Bbb Z$ will be called an interval if it is of the
form $I_{p,q}$. 

\vs {\bf Theorem 2.2.} {\it Let $\nu\in R$ and assume
$V$ is a finite-dimensional $\g$-module with respect to a
representation $\pi$. Let $\gamma$ be a $\tt$-weight of $V$ and let $$I
= \{j\in \Bbb Z\mid \gamma + j\,\nu\,\,\hbox{be a $\tt$-weight of
$V$}\}, $$ noting that $I$ is of course finite and not empty since $0\in
I$. Then there exist $p\leq 0\leq q,\,\,p,q\in \Bbb Z$ such that $$I =
I_{p,q}.\eqno (2.2)$$ Moreover if $I$ has only one element (i.e., $p=q=0$),
then $(\gamma,\nu)=0$. Furthermore if $I$ has more than one element
(i.e.,
$p<q$),  then $$\eqalign{(\gamma +
q\,\nu,\nu)&>0\,\,\hbox{and}\cr  (\gamma + p\,\nu,\nu)&<0.\cr}\eqno
(2.3)$$ Finally let $m\in I_{p,q}$. If $m<q$. Then
$$\pi(\g_{\nu})(V_{\gamma + m\,\nu})\neq 0, \eqno (2.4)$$ and if $p<m$, 
then $$\pi(\g_{-\nu})(V_{\gamma + m\,\nu})\neq 0.\eqno (2.5)$$}

{\bf
Proof.} Let $X = \sum_{j\in I}\,V_{\gamma + j\,\nu}$ so that $X$ is
stable under $\pi(\g_{\nu})$ and $\pi(\g_{-\nu})$ as well as
$\pi(\m)$.  One notes that, by (1.30), if $j\in I$, then $$V_{\gamma +
j\,\nu}\,\,\hbox{is the eigenspace of $\pi(h_{\nu})|X$ corresponding
to the eigenvalue $\langle \gamma,h_{\nu}\rangle + 2j$}.\eqno (2.6)$$
Now since $h_{\nu} \in [\g_{\nu},\g_{-\nu}]$ one must have
$$\hbox{tr}\,\pi(h_{\nu})|Y = 0\eqno (2.7)$$ where $Y \s X$ is any subspace
which is stable under $ \pi(\g_{\nu})$ and $\pi(\g_{-\nu})$. If $I$
has one element, then obviously $I= I_{0,0}$ and one has 
$(\gamma,\nu)=0$ by (2.7). Thus it suffices to
consider the case where
$I$ has more than one element. Now if $Y_1,Y_2$ are two nonzero
subspaces of $X$ that are both stable under $ \pi(\g_{\nu})$ and
$\pi(\g_{-\nu})$, it follows from (2.7) that one cannot have that the 
$$\hbox{
maximal eigenvalue of $\pi(h_{\nu})$ in $Y_1$ < minimal
eigenvalue of $\pi(h_{\nu})$ in $Y_2$}.\eqno (2.8)$$ Now assume that
$p,q\in I$ and $m\in \Bbb Z$ is such that $m\notin I$ and $p< m<q$. If
we define $Y_1$ (resp. $Y_2$) to be the sum of all $V_{\gamma +
j\,\nu}$,  where $j\in I$ and $j<m$ (resp. $j>m$),  the conditions of
(2.8) are satisfied which, as noted above, is a contradiction. Thus $I=
I_{p,q}$ for some $p,q\in \Bbb Z$ where 
$q> p$. But then (2.3) follows from (2.7) where $Y= X$.

 Now let $m\in I_{p,q}$ where $m<q$. Assume that
$\pi(\g_{\tau})(V_{\gamma + m\,\nu}) = 0$. That is, 
$$\pi(e_{\varphi})(V_{\gamma + m\,\nu})= 0,\,\,\forall \varphi\in
\Delta(\g_{\nu}).\eqno (2.9)$$ Thus for any $v\in V_{\gamma +
(m+1)\,\nu}$ and $\varphi \in \Delta(\g_{\nu})$, one has
$$\pi(e_{\varphi})\pi(e_{-\varphi})v=0.\eqno (2.10)$$
But from the representation theory of three-dimensional simple Lie algebras, (2.10) implies that
$\pi(e_{-\varphi})v=0$. That is $$\pi(\g_{-\nu})(V_{\gamma +
(m+1)\,\nu})= 0. \eqno (2.11)$$ But then if $Y_1$ (resp. $Y_2$) is the
sum of all $V_{\gamma + j\,\nu}$ for $j\in I_{p,q}$ where $j\leq m$
(resp. $j\geq m+1$), one defines $Y_1$ and $Y_2$ satisfying the
contradictory (2.8). This proves (2.4). Clearly a similar argument
proves (2.5). QED

 \vs  Applying Theorem 2.2 and Lemma 2.1 to the case where $V=\g$ and
$\pi$ is the adjoint representation, one immediately has the following result
asserting that some familar properties of ordinary roots still hold for $\tt$-roots.
 
\vs {\bf Theorem 2.3.} {\it Let $\nu,\mu \in R$. If $\mu +\nu\in R$ (resp. $\mu-\nu\in R$), then 
$$\eqalign{[\g_{\mu},\g_{\nu}] &= \g_{\nu +\mu}\,\,(\hbox{resp.}\cr [\g_{\mu},\g_{-\nu}] &= \g_{\mu
-\nu}).\cr}\eqno (2.12)$$ Furthermore, one indeed has $$\mu + \nu\in R\,\, (\hbox{resp.}\,\,\mu - \nu\in R)\,\,
\hbox{if}\,\,
(\mu,\nu)<0\,\, (\hbox{resp.}\,\,
(\mu,\nu)>0)\,\,\hbox{and}\,\, \mu +\nu \neq 0 \,(\hbox{resp.}\,\, \mu-\nu\neq 0).\eqno (2.13)$$ Moreover 
$$\hbox{if}\,\,
(\mu,\nu) = 0\,\,\hbox{then}\,\, 
\mu +\nu\in R\,\,\hbox{if and only if}\,\,\mu-\nu\in R. \eqno (2.14)$$}
\vs
{\bf 2.2.} Recalling (1.3)
let $\delta_{\n}$ be in the dual space $\m^*$ to $\m$ defined so that if $x\in \m$, 
then $$\langle \delta_{\n},x\rangle = \hbox{\rm tr}\, \hbox{\rm ad}\,x|\n. \eqno (2.15)$$ Since $\n$ and
$\overline {\n}$ (see (1.10)) are clearly stable under $\hbox{\rm ad}\,\m$ one has a partition $R
= R_{\n}\cup R_{\overline {\n}}$ so that $$\eqalign{\n&= \sum_{\nu\in
R_{\n}}\g_{\nu}\cr \overline{\n}&= \sum_{\nu\in R_{\overline {\n}}}\g_{\nu}.\cr}\eqno
(2.16)$$ Clearly Remark 1.4 implies $$R_{\overline {\n}} = - R_{\n}.\eqno (2.17)$$

 {\bf Lemma 2.4.} {\it One has $\delta_{\n}\in \tt_{\Bbb R}^*$ (see (1.7)).
Furthermore
$$\eqalign{(\delta_{\n},\nu) &>0,\,\,\hbox{if $\nu\in R_{\n}$ and}\cr
(\delta_{\n},\nu) &<0,\,\,\hbox{if $\nu\in R_{\overline{\n}}$}. \cr}\eqno (2.18)$$}

{\bf Proof.} Since $\ss$ is semisimple,  $\delta_{\n}$ must vanish on $\ss$ and
hence $$\delta_{\n} \in \tt^*\eqno (2.19)$$ where, besides regarding $\tt^*\s \hh^*$ as in
(1.6), we also regard
$\tt^*\s
\m^*$, using (1.4). Let $\varphi\in \Delta$. Normalize the choice of root vectors so
that $(e_{\varphi},e_{-\varphi}) = 1$. Let $\a_{\varphi}$ be the root TDS
corresponding to $\varphi$ and let $x_{\varphi} = [e_{\varphi},e_{-\varphi}]$. Then
as one knows $x_{\varphi}\in \hh^*_{\Bbb R}$ and (see (1.27)) $$\tau(x_{\varphi}) =
\varphi.\eqno (2.20)$$ For some index set $P$, let $$\g = \sum_{p\in P}\u_p\eqno (2.21)$$
be the decomposition of $\g$ into a sum of irreducible $\hbox{\rm ad}\,(\a_{\varphi} +
\hh)$-submodules. One clearly has the direct
sum (with possibly $0$-dimensional summands)
$$\n=
\sum_{p \in P} \n_p\eqno (2.22)$$ where $\n_p = \n\cap \u_p$. Now assume that
$\varphi\in \Delta(\n)$. Then since each
$\n_p$ is
stable under
$\hbox{\rm ad}\,e_{\varphi}$ for any $p\in P$, it is immediate from the representation theory of a
TDS that
$$\hbox{\rm tr}\,\hbox{\rm ad}\,x_{\varphi}|\n_p\geq 0.$$ But there exists $p_o\in P$ such that $\n_{p_o} = \Bbb
C e_{\varphi}$ so that
$\hbox{\rm tr}\, \hbox{\rm ad}\,x_{\varphi}|\n_{p_o}>0$. Thus $$\delta_{\n}(x_{\varphi})>0.\eqno (2.23)$$ 
But
then
$$(\delta_{\n},\varphi)>0\eqno (2.24)$$ by (2.20). Now let $\nu\in R_{\n}$ and let
$\varphi\in \Delta(\g_{\nu})$. But clearly $\nu = \varphi_{\tt}$ (see (1.23)) so that
(2.24) and (2.19) imply the first line of (2.18). The second line is implied by (2.17).
Since $\hh^*_{\Bbb R}$ is clearly spanned by $R_{\n}$ (see (1.28)), it follows from (2.18)
and (2.10) that $\delta_{\n} \in \tt^*_{\Bbb R}$. QED

\vs Introduce the lexicographical ordering
in $\tt_{\Bbb R}^*$ with respect to an orthogonal ordered basis of $\tt^*_{\Bbb R}$
having
$\delta_{\n}$ as its first element. It follows  from Lemma 2.4 that if $R_+$ is the
set of positive $\tt$-roots with respect to this ordering, one has $$R_+ = R_{\n}.\eqno
(2.25)$$
\vskip .5pc {\bf Remark 2.5.} One recalls that since $\tt^*_{\Bbb R}$ is a lexicographically ordered real
euclidean space if
$\xi_i\in \tt^*_{\Bbb R},\,i=1,\ldots,k$, are positive elements such that, for $i\neq j$,
$$(\xi_i,\xi_j)\leq  0,$$ then the $\xi_i$ are linearly independent. (
See e.g., [Hum], \S 10, Theorem$^{\prime}$, (3), p.~48.)

\vs {\bf 2.3.} Let
$\ell(\tt) = \hbox{\rm dim}\,\tt$ and $\ell(\ss) = \hbox{\rm dim}\,\hh(\ss)$ so that $\ell(\ss)$ is the rank of
$\ss$ and $$\ell = \ell(\tt) + \ell(\ss).\eqno (2.26)$$ Let $\Pi =
\{\alpha_1,\ldots,\alpha_{\ell}\}$ be the set of simple positive roots in $\Delta_+$. If
$\varphi\in \Delta_+(\ss)$ and $\varphi = \varphi_1 + \varphi_2$ where $\varphi_1,\,
\varphi_2\in \Delta_+$, obviously $(\varphi_1)_{\tt} = -(\varphi_2)_{\tt}$. But then
$(\varphi_i)_{\tt},\,\,i=1,2$, vanish by Lemma 2.4 so that $\varphi_1, \varphi_2\in
\Delta_+(\ss)$. Hence if $\varphi\in \Delta_+(\ss)$ is simple with respect to
$\Delta_+(\ss)$, it is simple with respect to $\Delta_+$. We may therefore order $\Pi$ so
that $\alpha_{\ell(\tt) +i} \in \Delta_+(\ss)$, for $i=1\,\ldots,\ell(\ss)$, and hence if
$\Pi_{\ss} =
\{\alpha_{\ell(\tt)+1},\ldots,\alpha_{\ell}\}$, then $$\Pi_{\ss}\,\,\hbox{is a basis of
$\hh(\ss)^*$; see (1.6).}\eqno (2.27)$$ \vskip .5pc A $\tt$-root $\nu\in
R_+$ is called simple if $\nu$ cannot be written $\nu = \nu_1 +\nu_2$ where
$\nu_1,\nu_2\in R_+$. Let $R_{\hbox{\smallrm simp}}$ be the set of simple $\tt$-roots in $R_+$.

 \vs
{\bf Lemma 2.6.} {\it Assume $\xi_1,\xi_2\in R_{simp}$ are distinct. Then
$$(\xi_1,\xi_2)\leq 0. \eqno (2.28)$$ so that, by Remark 2.5, the elements in $R_{\hbox{\smallrm simp}}$
are linearly independent. In particular $$\hbox{\rm card}\,R_{\hbox{\smallrm simp}} \leq \ell(\tt).\eqno
(2.29)$$}

 \vskip .3pc
{\bf Proof}. Assume $(\xi_1,\xi_2)>0$. Then, by Theorem 2.3, $\xi_1-\xi_2$ and $\xi_2-\xi_1$
are in $R$. Without loss assume $\nu\in R_+$ where $\nu = \xi_1-\xi_2$. Then $\xi_1 = \nu
+ \xi_2$. This contradicts the simplicity of $\xi_1$. Hence one has (2.28). QED

\vs As
noted in (1.6), $\hh(\ss)^*$ is the orthocomplement of $\tt$ in $\hh$ . Thus, by 
(2.27), if $$\beta_j = (\alpha_{j})_{\tt} \eqno (2.30)$$ for
$j=1,\ldots,\ell(\tt)$, then clearly  $$\eqalign{\beta_j& \in R_+\,\,\hbox{and}\cr
\beta_j&,\,\, j=1,\ldots,\ell(\tt),\,\,\hbox{are a basis of $\tt^*_{\Bbb R }$}.\cr}\eqno
(2.31)$$ 

 Recalling (1.1), for any $\nu\in R$ and $j=1,\ldots, \ell(\tt)$,  let
$n_j(\nu) = n_{\alpha_j}(\varphi)$ where $\varphi\in \Delta(\g_{\nu})$. This is
independent of the choice of $\varphi$ by (2.27). \vs {\bf Theorem 2.7.} {\it One has
$$R_{\hbox{\smallrm simp}} =
\{\beta_1,\ldots,\beta_{\ell(\tt)}\}\eqno (2.32)$$ so that $R_{\hbox{\smallrm simp}}$ is a basis of
$\tt^*_{\Bbb R}$ and for
$i\neq j$,
$$(\beta_i,\beta_j)\leq 0.\eqno (2.33)$$ Furthermore for $\nu\in R_+$ and
$j=1,\ldots,\ell(\tt)$, one has $$\nu =
\sum_{j=1}^{\ell(\tt)} \,n_j(\nu)\,\beta_j.\eqno (2.34)$$}

 {\bf Proof.} Let $j\in
\{1,\,\ldots,\ell(\tt)\}$. Assume that $\beta_j\notin R_{\hbox{\smallrm simp}}$. Then there exists
$\nu_1,\nu_2\in R_+$ such that $\beta_j = \nu_1 +\nu_2$. But then by (3) in Theorem 2.3
one has $$[\g_{\nu_1},\g_{\nu_2}] = \g_{\beta_j}.\eqno (2.35)$$ But $e_{\alpha_{j}}\in
\g_{\beta_j}$ by (1.23). Thus there exist, for $i=1,2$, $\varphi_i\in
\Delta(\g_{\nu_i})$ such that $\varphi_1 + \varphi_2 = \alpha_{j}$. This
contradicts the simplicity of $\alpha_{j}$ since $\varphi_1 , \varphi_2\in
\Delta_+$. But then (2.32) follows from (2.29). Also (2.33) and the fact that $R_{\hbox{\smallrm simp}}$
is a basis of $\tt^*_{\Bbb R}$ follow from Lemma 2.6. 

Now let $\nu\in R_+$ and let $\varphi\in \Delta(\g_{\nu})$. Recalling the expansion (1.1)
one immediately has $$\nu = \sum_{\alpha\in \Pi} \,n_{\alpha}(\varphi)\,\alpha_{\tt},\eqno
(2.36)$$ but this yields (2.34). QED

\vs {\bf Remark 2.8.} Note that if $j=1,\ldots,
\ell(\tt)$, one has that 
$e_{\alpha_{j}}\in
\g_{\beta_j}$ and $e_{\alpha_{j}}$ is a lowest weight vector for the
irreducible $\m$-module $\g_{\beta_j}$. Indeed this is clear since for
$i=1,\ldots,\ell(\ss)$,
$$[e_{-\alpha_{\ell(\tt) +1}},e_{\alpha_j}]= 0.\eqno (2.37)$$

\vskip .5pc {\bf 2.4.} Henceforth we will assume that $\g$ is simple. Let $\psi\in \Delta_+$
be the highest root so that, by the simplicity of
$\g$, 
$$\hbox{\rm cent}\,\n_{\b} = \Bbb C\,e_{\psi} \eqno (2.38)$$ where we recall (\S 1.1) $\n_{\b} =
[\b,\b]$. \vskip .5pc {\bf Remark 2.9.} For any
$\alpha\in \Pi$ one knows that $$n_{\alpha}(\psi)>0\eqno (2.39)$$ and, for any $\varphi\in
\Delta_+$, $$n_{\alpha}(\psi)\geq n_{\alpha}(\varphi).\eqno (2.40)$$ Indeed (2.39) and
(2.40) are consequences of the immediate fact that $$e_{\psi}\in
U(\n_{\b})\,e_{\varphi}\eqno(2.41)$$ where any Lie algebra $\a$, $U(\a)$ is the enveloping
algebra of $\a$.\vs 
For any
$\nu\in R$ let
$o(\nu) =
\sum_{j=1}^{\ell(\tt)}\, n_{\beta_j}(\nu)$, and for any $k\in \Bbb Z_+$ let $$\n(k) =
\sum_{\nu\in R,\,o(\nu)= k}\,\g_{\nu}.\eqno (2.42)$$ For $j,k\in \Bbb Z$,  clearly
$$[\n(j),\n(k)] \s \n(j+k).\eqno (2.43) $$ Now $\psi_{\tt} \neq 0$ by Remark 2.9 (one
has $\ell(\tt) < \ell$ by our assumption in \S 1.1). Let
$\nu(\hbox{\rm cent})\in R$ be defined by putting $\nu(\hbox{\rm cent}) = \psi_{\tt}$. Put $k(\hbox{\rm cent})=
o(\nu(\hbox{\rm cent}))$ so that
$$k(\hbox{\rm cent}) = \sum_{i=1}^{\ell(\tt)} n_{\alpha_i}(\psi).\eqno (2.44)$$ Remark 2.9 clearly also
implies
\vs {\bf Proposition 2.10.} {\it $\n(k) = 0 $ if $k> k(\hbox{\rm cent})$ and $$\n(k(\hbox{\rm cent})) =
\g_{\nu(\hbox{\rm cent})}.\eqno (2.45)$$ Furthermore one has the direct sum $$\n = \sum_{k=0}^
{k(\hbox{\smallrm cent})} \n(k).\eqno (2.46)$$}

 The upper central series of $\n$ (defined for any
nilpotent Lie algebra) is a sequence of distinct ideals $\n_1\s n_2 \s
\cdots
\s
\n_d =
\n$ where
$\n_1 = \hbox{\rm cent}\,\n$ and for $i\geq 2$, $$n_i/\n_{i-1} = \hbox{\rm cent}\,\n/\n_{i-1}.\eqno (2.47)$$ See
(14) on p. 29 in [Jac].
We refer to $d$ as the length of the upper central series.\vs {\bf Theorem 2.11.} {\it If $i<
k(\hbox{\rm cent})$ and
$\nu\in
\n(i)$, there exists
$j\in
\{1,\ldots,\ell(\tt)\}$ such that $\nu + \beta_j\in \n(i+1)$ so that 
$$[\g_{\nu},\g_{\beta_j}]= \g_{\nu + \beta_j}.\eqno (2.48)$$ In particular
$[\g_{\nu},\g_{\beta_j}]\neq 0$. Furthermore $$\g_{\nu(\hbox{\rm cent})} = \hbox{\rm cent}\,\n\eqno (2.49)$$ and
the upper central series $\n_i$, of the nilradical $\n$ of the general proper parabolic
subalgebra
$\q$, is given as follows: $$\n_i = \sum_{j=1}^i\n(k(\hbox{\rm cent}) -j +1)\eqno(2.50)$$ noting that
$k(\hbox{\rm cent})$ is the length of the upper central series.}

\vs {\bf Proof.} Let $\varphi\in
\Delta(\g_{\nu})$ be the highest weight of the $\m$-irreducible module $\g_{\nu}$. Thus
$[e_{\alpha_{i+\ell(\tt)}},e_{\varphi}]=0$ for $i=1,\ldots,\ell(\ss)$. But $\varphi\neq
\psi$. Hence by the uniquenss of the $\g$-highest weight there must exist
$j=1,\ldots,\ell(\tt)$ such that $[e_{\alpha_j},e_{\varphi}]\neq 0$. But this implies
(2.48). Now obviously $\g_{\nu(\hbox{\rm cent})}\s \hbox{\rm cent}\,\n$ by (2.43) and (2.46). But clearly
$\hbox{\rm cent} \,\n$ is stable under $\hbox{\rm ad}\,\m$. Therefore to prove (2.49) it suffices to show that
$\g_{\nu}\not\subset \hbox{\rm cent}\,\n$ for $\nu\neq \nu(\hbox{\rm cent})$. But this is established by (2.48).
Let $i>1$. Assume inductively that one has (2.50) where $i-1$ replaces $i$. Then, returning
to (2.50) as stated, if
$\n_{(i)}$ is given by the right side of (2.50),  one has $\n_{(i)}\s
\n_{i}$ by (2.43) and (2.46). But the upper central series (e.g.,  by induction) is
stabilized by $\hbox{\rm ad}\,\m$ so that in particular $\n_i$ is stabilized by $\hbox{\rm ad}\,\m$. But again
(2.48) implies that if $\g_{\nu}\s \n(j)$ where $j< k(\hbox{\rm cent}) -i +1$, then $\g_{\nu}
\not\subset \n_i$. This implies that
$\n_{i} = \n_{(i)}$. QED

\vs  Using the notation in [Jac], see page 23, the lower
central series $\n^i$ of $\n$ is a sequence of ideals defined inductively so that $\n^1 =
\n$ and for $i>1, \n^i = [\n,\n^{i-1}]$. See also p.~11 in [Hum]. The indexing in [Hum] differs by 1 from
the indexing in [Jac]. We will call the maximum $k$ such that $\n^k\neq 0$ the length of the lower central
series.

 \vs {\bf Theorem 2.12. } {\it Let $i$ be any integer where $2\leq i\leq k(\hbox{\rm cent})$ and let $\nu
\in R$ where $\g_{\nu}\in \n(i)$. Then there exists $j\in \{1,\ldots,\ell(\tt)\}$ and
$\mu\in R$ where $\g_{\mu}\in \n(i-1)$ such that $$[\g_{\beta_j},\g_{\mu}] = \g_{\nu}\eqno
(2.51)$$ so that $$[\n(1),\n(i-1)] = \n(i).\eqno (2.52)$$ In particular the lower central
series
$\n^i$ of the nilradical $\n$ of the arbitrary proper parabolic subalgebra $\q$ of $\g$ is
given by $$\n^i = \sum_{j=i}^{k(\hbox{\smallrm cent})} n(j)\eqno (2.52)$$ so that (see 2.44) $k(\hbox{\rm cent})$
is the length of the lower central series of $\n$ (as well as the upper central series, see
Theorem 2.11). The lower and upper central series of $\n$ are therefore, except for
indexing, the same $$\n^i = \n_{k(\hbox{\smallrm cent})-i +1}\eqno (2.53)$$ for
$i=1,\ldots,k(\hbox{\rm cent})$}.

\vs
 {\bf Proof.} It suffices only to prove (2.51). But, by Remark 2.5
and Theorem 2.7, $\beta_j\in R_{\hbox{\smallrm simp}}$ exists so that $$(\nu,\beta_j) > 0.\eqno (2.54)$$
Thus $\mu\in R$ where $\mu = \nu-\beta_j$ by Theorem 2.3. But also
necessarily $n_{j}(\nu)>0$ by (2.33) so that $\g_{\mu}\in \n(i-1)$. But $\mu + \beta_j =
\nu$. Thus by (3) of Theorem 2.3 one has (2.51). QED

\vskip 1.5pc \centerline{\bf 3. Borel--de Siebenthal theory, special elements,}

\centerline{\bf  and the Lie subalgebras they define}\vskip 1pc

 {\bf 3.1.} We continue to assume (starting in \S
2.4) that $\g$ is simple. In this section we will apply the results of \S 1 and \S 2 to the case where
$\ell(\tt) = 1$. It will be convenient to change some notation and earlier indexing. In particular we now fix
an ordering in $\Pi$ so that $\Pi = \{\alpha_1,\ldots,\alpha_{\ell}\}$. Also recalling (1.1)
we will write $n_{i}(\varphi)$ for $n_{\alpha_i}(\varphi)$ so that for the highest root
$\psi$ one has $$\psi = \sum_{i=1}^{\ell} \,n_i(\psi)\,\alpha_i.\eqno (3.1)$$ \vskip .5 pc
Now let $x_j,\,j=1,\ldots,\ell$, be the basis of $\hh_{\Bbb R}$ so that
$$\langle\alpha_l,x_j\rangle = \delta_{ij}.\eqno (3.2)$$ Now for $j=1,\ldots,\ell$, let
$\q[j]$ be the maximal standard parabolic subalgebra defined by $x_j$. The standard Levi
factor $\m[j]$ of $\q[j]$ is the centralizer $\g^{x_j}$ of $x_j$ in $\g$. The decomposition
(1.3) now becomes $$\m[j] = \tt[j] + \ss[j]\eqno (3.3)$$ where the one-dimensional
center $\tt[j]$ of $\m[j]$ is given by $$ \tt[j] = \Bbb C\,x_j.\eqno (3.4)$$
If $\Pi_{\ss[j]}$ is the set of simple positive roots of the (rank $\ell-1$)-semisimple
Lie algebra $\ss[j]$, defined as in (2.27), one now has $$\Pi_{\ss[j]} =
\Pi\setminus \{\alpha_j\}.\eqno (3.5)$$ The nilradical $\n[j]$ of $\q[j]$ is given by $$\n[j]
= \hbox{span of}\,\,\{e_{\varphi}\mid n_j(\varphi) > 0\}. \eqno (3.6)$$

 Write $R[j]$ 
for $R$, $R[j]_+$ for  $R_+$ and $R[j]_{\hbox{\smallrm simp}} $ for $R_{\hbox{\smallrm simp}}$. Let $\beta[j] =
(\alpha_j)_{\tt[j]}$ so that $$\langle \beta[j],x_j\rangle = 1.\eqno (3.7)$$\vskip .5pc {\bf Proposition 3.1.}
{\it One has
$$R[j]_{\hbox{\smallrm simp}} = \{\beta[j]\}\eqno (3.8)$$ and $$R[j]_+ = \{\beta[j],
2\,\beta[j],\ldots,n_j(\psi)\,\beta[j]\}\eqno (3.9) $$ and $$R[j] = \{\pm\beta[j],
\pm 2\,\beta[j],\ldots,\pm n_j(\psi)\,\beta[j]\}.\eqno (3.10) $$ In particular
$$\eqalign{\hbox{\rm card}\,R[j]_+  &= n_j(\psi)\cr \hbox{\rm card}\,R[j] 
&= 2\,n_j(\psi).\cr}\eqno (3.11)$$}\vs {\bf Proof. } The proof is immediate from (2.17), (2.25), 
Theorem 2.7, Theorem 2.12 and (2.44) which implies here that $$k(\hbox{\rm cent}) = n_j(\psi).\eqno (3.12)$$
\hfill QED\vs

 Let 
$I[j] = \{\pm 1,\ldots,\pm n_j(\psi)\}$ so that if $\nu\in R[j]$, then $\nu = k\,\beta[j]$ for $k\in
I[j]$. We recall (see \S 1.2) that for $k\in I[j]$, $$\eqalign{\g_{k\beta[j]}&= \hbox{\rm span
of}\,\,\{e_{\varphi}\mid n_j(\varphi) = k\}\cr &= \{z\in \g\mid [x_j,z] = k\,z\}.\cr}\eqno (3.13)$$
One has the direct sums $$\eqalign{\n[j] &= \sum_{k=1}^{n_j(\psi)}\,\g_{k\beta[j]}\cr \g &=\m[j]
+ \sum_{k\in I[j]} \,\g_{k\,\beta[j]}.\cr}\eqno (3.14) $$\vskip .5pc One now has \vs 

{\bf Theorem 3.2.}
{\it Let $j\in \{1,\ldots,\ell\}$. Then $\g_{k\beta[j]}$ is $\hbox{\rm ad}\,\m[j]$-irreducible for any $k\in
I[j]$. In particular the nilradical $\n[j]$ of $\q[j]$ is, as an $\hbox{\rm ad}\,\m[j]$-module, 
 multiplicity-free with $n_j(\psi)$-irreducible components. Also $\g/\m[j]$ is, as an $\hbox{\rm
ad}\,\m[j]$-module, multiplicity-free with $2\,n_j(\psi)$-irreducible components. Finally, if $p,q\in I[j]$ and
$p +q \in I[j]$,  then $$[\g_{p\beta[j]}, \g_{q\beta[j]}] = \g_{(p+q)\beta[j]}.\eqno
(3.15)$$}

 {\bf Proof.} The first statement is just the present application of Theorem 1.9. The
equality (3.15) is given by (3) of Theorem 2.3. QED

\vs {\bf 3.2.} 
Let $C\s\hh_{\Bbb R}$ be the fundamental Weyl chamber corresponding to $\b$ so that
$$C=\{x\in \hh_{\Bbb R}\mid \langle \alpha_i,x\rangle \geq
0,\,\,\hbox{for}\,\,i=1,\ldots,\ell\},\eqno (3.16)$$ and let $A\s C$ be the fundamental
alcove so that $A$ is the simplex defined by $$A=\{x\in C\mid \langle \psi,x\rangle\leq
1\}.\eqno (3.17)$$ 

Let $G$ be a simply-connected complex group for which $\g =\hbox{\rm Lie}\,G$. Let $K$ be a maximal compact
  subgroup  of $G$. We may choose $K$ so that if $\k = \hbox{\rm Lie}\, K$,  then $i\,h_{\Bbb R}=$ 
$\hbox{\rm Lie}\,H_{\hbox{\smallrm comp}}$ where $H_{\hbox{\smallrm comp}}$ is a maximal torus of $K$. 

A classical result of Cartan and Weyl is the statement \vs 

{\bf
Proposition 3.3.} {\it For any element $g\in K$ there exists a unique element $x\in A$ such
that $$g\,\, \hbox{is $K$-conjugate to $\hbox{\rm exp}\,2\,\pi i \ x$}. \eqno (3.18)$$}\vs

 Clearly the
$\ell+1$ vertices $v_j,\,j=0,\ldots,\ell$, of $A$ are then given as $v_0 = 0$ and for
$j>0$, $$v_j= x_j / n_{j}(\psi).\eqno (3.19)$$ For $j=1,\ldots, \ell$, let $a_j\in K$ be defined by putting
$$a_j =\hbox{\rm exp}\,2\,\pi \ i\,v_j.\eqno (3.20) $$ Let $\omega[j]$ be the $n_j(\psi)$-root of unity given by
putting $\omega[j]= e^{2\,\pi\,i/n_j(\psi)}$. For $g\in G$ let $G^g$ (resp. $\g^g$) be the centralizer of
$g$ in $G$ (resp. $\g$). As an immediate consequence of (3.14) the well-known adjoint action of
$a_j$ on $\g$ is given by
 
\vs {\bf Proposition 3.4.} {\it $\hbox{\rm Ad}\,a_j$ has order $n_j(\psi)$ on $\g$. In fact one has
$$\eqalign{\hbox{\rm Ad}\,a_j &= 1\,\,\hbox{on}\,\,\m[j]\cr &=
\omega[j]^k\,\,\hbox{on $\g_{k\,\beta[j]}$ for all $k\in I[j]$}.\cr}\eqno (3.21)$$ In particular
one has the direct sum $$\g^{a_j} = \m[j] + \g_{n_j(\psi)\,\beta[j]} + \g_{-n_j(\psi)\,\beta[j]}.\eqno
(3.22)$$}\vs

 An element $a\in K$ is called special if the reductive subalgebra $\g^a$ of $\g$ is in fact
semisimple (i.e.,  $\hbox{\rm cent}\,\g^a= 0$). The following is also well known but proved here for
completeness.\vs 

{\bf  Proposition 3.5.} {\it For $j= 1,\ldots, \ell$, the element $a_j$ is special.}\vs {\bf
Proof.} Since
$\hh\s \m[j]$ one obviously has $\hbox{\rm cent}\,\g^{a_j} \s \hbox{\rm cent}\,\m[j]$. But
$\hbox{\rm cent}\,\m[j] = \Bbb C\,x_j$. However
$\hbox{\rm ad}\,x_j = n_j(\psi)$ on $\g_{n_j(\psi)\,\beta[j]}$. Thus $\hbox{\rm cent}\,\g^{a_j} = 0$. QED

\vs {\bf
Remark 3.6.} One can readily prove $a\in K$ is special if and only if either
$a=1$ or there exists $j\in \{1,\ldots,\ell\}$, necessarily unique, such that $a$ is $K$-conjugate to
$a_j$. We also remark that special elements arise in connection with distinguished nilpotent conjugacy
classes in 
$\g$. Indeed if $e$ is a distinguished nilpotent element then, where $G^e_o$ is the identity component of
$G^e$, the component group $G^e/G^e_o$ is isomorphic to the finite group $F$ where $F$ is the
centralizer in $G$ of a TDS containing $e$. Since $F$ is finite we may make choices so that $F\s K$. But
then the elements of $F$ are special. Indeed if $a\in F$,  then $e\in \g^a$. But if $x\in \hbox{\rm
cent}\,\g^a$,  then 
$x$ is a semisimple element that commutes with $e$. Thus $x=0$ since $e$ is distinguished. Hence
$\hbox{\rm cent}\,\g^a = 0$.

\vs {\bf 3.3.} Let $j= 1,\ldots,\ell$. For completeness in this section we wish to give the proof of the
 Borel--de Siebenthal algorithm for determining the Dynkin diagram of $\g^{a_j}$. We recall that the
extended Dynkin diagram of $\g$ is the usual Dynkin diagram (whose nodes are identified with $\Pi$) of $\g$,
together with an additional node $\alpha_0$, where if $i=1,\ldots,\ell$, then $\alpha_0$ is linked to
$\alpha_i$ with $$m_i = 2\,(\alpha_i,\psi)/(\alpha_i,\alpha_i)\eqno (3.23)$$ lines and arrowhead at
$\alpha_i$, directed at $\alpha_i$, if $m_i >1$. \vs

 {\bf Theorem 3.7.} (Borel--de Siebenthal algorithm). {\it The Dynkin diagram of $\g^{a_j}$ is the sub\-diagram
of the extended Dynkin diagram remaining after deleting
$\alpha_j$ (and all lines linked to $\alpha_j$) from the extended Dynkin diagram. }\vs

 {\bf Proof.} Let $\m[j]_+ =
\m[j]\cap \b$ so that $\m[j]_+$ is a Borel subalgebra, containing $\hh$, of $\m[j]$. Now
$\g_{-n_j(\psi)\beta[j]}$ is clearly a commutative nilpotent subalgebra of $\g^{a_j}$ which is stable (and
in fact irreducible) under $\hbox{\rm ad}\,\m[j]$. But then of course $\g_{-n_j(\psi)\beta[j]}$ is stable under
$\hbox{\rm ad}\,\m[j]_+$. Thus $$\b[j] = \m[j]_+ + \g_{-n_j(\psi)\beta[j]} \eqno (3.24)$$ is a solvable Lie
subalgebra of $\g^{a_j}$. But, by dimension, $\b[j]$ is then a Borel subalgebra of $\g^{a_j}$. Regarding 
$\Delta(\b[j])$ as a system of positive roots for $\g^{a_j}$, let $\Pi[j]$ be the corresponding set of
simple positive roots. But now $-\psi\in \Delta(\g_{-n_j(\psi)\beta[j]})$. (In fact
 clearly $e_{-\psi}$ is a lowest weight vector for the irreducible action of $\hbox{\rm ad}\,\m[j]$ on
$\g_{-n_j(\psi)\beta[j]}$.) But one notes it is immediate that $-\psi\in \Pi[j]$. Since $\g^{a_j}$ has
rank $\ell$ one has, recalling (3.5), $$\Pi[j] = \Pi_{\ss[j]} \cup \{-\psi\}.\eqno (3.25)$$ But then the
Dynkin diagram of $\g^{a_j}$ is the Dynkin diagram of $\ss[j]$, together with the node defined by $-\psi$,
where if $\alpha_i\in \Pi_{\ss[j]}$, then $-\psi$ (a long root) is linked to $\alpha_i$ by (see (3.23))
 $m_i$ lines and arrowhead at
$\alpha_i$, directed at $\alpha_i$, if $m_i >1$. But this is the Borel--de Siebenthal algorithm. QED

\vs {\bf 3.4.} Let $j = 1,\ldots, \ell,$ and let $I_o[j] = I[j]\setminus \{\pm n_j(\psi)\}$. Let $$\r[j] =
\sum_{k\in I_o[j]}\,\g_{k\,\beta[j]}\eqno (3.26)$$ so that by (3.14) and (3.22) one has the $B$-orthogonal
direct sum $$\g = \g^{a_j} + \r[j]. \eqno (3.27)$$ For $k=1,\ldots,n_j(\psi)-1$, let $$\g[a_j]^k = \{ z\in
\g\mid \hbox{\rm Ad}\,a_j(z) = \omega[j]^k\,z\}\eqno (3.28)$$ Obviously $$\g[a_j]^k\,\,\hbox{is stable under
$\hbox{\rm ad}\,\,\g^{a_j}$ for $k=1,\ldots,n_j(\psi)-1$}.\eqno (3.29)$$ But clearly, by Proposition 3.4, for
$k=1,\ldots,n_j(\psi)-1$, 
$$\g[a_j]^k =
\g_{k\,\beta[j]} + \g_{(k-n_j(\psi))\,\beta[j]}, \eqno (3.30)$$ and hence by (3.26) one has the direct
sum $$\r[j] = \sum_{k=1}^{n_j(\psi)-1} \g[a_j]^k. \eqno (3.31)$$ Consequently one notes that not only is
$\r[j]$ stable under $\hbox{\rm ad}\,g^{a_j}$ but (3.31) isolates $n_j(\psi)-1$
$\hbox{\rm ad}\,g^{a_j}$-submodules of
$\r[j]$.
\vs 

{\bf Theorem 3.8.} {\it Let $j=1,\ldots, j$, so that $a_j$, defined by the vertex $v_j$ of the
fundamental alcove, is a special element of $K$. In particular,  its centralizer $\g^{a_j}$ is a
maximal semisimple Lie subalgebra of $\g$ if $n_j(\psi)$ is prime---by Borel--de Siebenthal theory. Let
$\r[j]$ be the
$B$-orthocomplement of $\g^{a_j}$ in $\g$. Then $\r[j]$ is a multiplicity-free $\hbox{\rm ad}\,\m[j]$-module
where
$\m[j]\s \g^{a_j}$ is the centralizer of $v_j$ in $\g$ and (3.26) is the decomposition of $r[j]$ into a
sum of $2(n_j(\psi)-1)$-irreducible $\hbox{\rm ad}\,\m[j]$-submodules. Next (3.30) defines the decomposition of
the $\hbox{\rm Ad}\,a_j$-weight space $\g[a_j]^k$ into a sum of two irreducible $\hbox{\rm
ad}\,m[j]$-submodules. Furthermore, and mainly, the $\hbox{\rm Ad}\,a_j$-weight space $\g[a_j]^k$ is an
irreducible $\hbox{\rm ad}\,\g^{a_j}$-submodule. In addition  $\r[j]$ is a multiplicity-free
$\hbox{\rm Ad}\,\g^{a_j}$-module and (3.31) is the decomposition of
$r[j]$ into a sum of $n_j(\psi)-1$-irreducible $\hbox{\rm ad}\,\g^{a_j}$-submodules. Finally, if $p,q,r\in
\{1,\ldots,n_j(\psi)-1\}$ and $r\cong p+q\,\, \hbox{mod}\,\,n_j(\psi)$,  then $$ [\g[a_j]^p,\g[a_j]^q]=
\g[a_j]^r.\eqno (3.32)$$}

 {\bf Proof.} Up until the sentence beginnning with ``Furthermore" the stated
results have been established in Theorem 3.2. But now
$$[\g_{(k-n_j(\psi))\,\beta[j]},\g_{n_j(\psi)\,\beta[j]}] = \g_{k\,\beta[j]}$$ and
$$[\g_{k\,\beta[j]},\g_{-n_j(\psi)\,\beta[j]}]= \g_{(k-n_j(\psi))\beta[j]}$$ by (3.15). Hence
recalling (3.22) and (3.30) it follows that $\g[a_j]^k$ is $\g^{a_j}$-irreducible for all
$k=1,\ldots,n_j(\psi)-1$. 

Clearly the left side of (3.32) is contained in the right side of (3.32). But by $\hbox{\rm ad}\,g^{a_j}$
irreducibility one has (3.32) as soon as one observes that the left side is nonzero. But this is clear
from (3.15) if $p+q = r$. If $p + q > n_j(\psi)$, then $r = p+q -n_j(\psi)$. But then (3.32) follows from
(3.15) where $q$ is replaced by $q-n_j(\psi)$. QED

\vs {\bf Remark 3.9.} One of the main results of  Borel--de Siebenthal theory is the statement that $\g_1$ is
a maximal proper (i.e., $g_1\neq \g$) semisimple subalgebra such that
$\hbox{\rm rank}\, \g_1 = \hbox{\rm rank}\,\g$ if and only if $$\g_1 \cong \g^{a_j}\eqno (3.33)$$ where
$n_j(\psi)$ is a prime number. This may be proved as follows: In one direction assume (3.33) where $n_j(\psi)$
is prime. Then if $\g_o$ is a subalgebra where
$\g^{a_j}\s \g_o$ and $\g_o\neq \g^{a_j}$,  there must exist, by (3.27), $0\neq x\in \r[j]\cap \g_o$. But since
$\g_o$ is stable under $\hbox{\rm ad}\,\g^{a_j}$ it follows from Theorem 3.8 that there exists 
$k\in \{1,\ldots, n_j(\psi)-1\}$
such that $\g[a_j]^k\s \g_o$. But $\g$ is generated by $\g^{a_j}$ and $\g[a_j]^k$, by (3.32), since $n_j(\psi)$
is prime. Thus $\g^{a_j}$ is maximal as a proper Lie subalgebra of $\g$. Conversely assume $\g_1$ 
is a maximal proper semisimple subalgebra where
$\hbox{\rm rank}\,\g_1 = \hbox{\rm rank}\,\g$. Let $G_1\s G$ be the subgroup corresponding to $\g_1$. Let
$\gamma$ denote the adjoint  representation of
$G_1$ on $\g/\g_1$. By the equal rank condition $0$ is not a weight of $\gamma$. Thus $\gamma$ does not descend
to the adjoint group of $\g_1$. Thus there exists $1\neq c\in \hbox{\rm cent}\,G_1$ such that $c\notin
 \hbox{\rm Ker}\,\gamma$. But $c$ has finite order since $\g_1$ is semisimple. We may therefore  make choices so
that $c\in K$. Of course
$\g_1\s \g^c$. By maximality
$$\g^c = \g_1.\eqno (3.34)$$ But then $c$ is special and by Remark 3.6 choices can be made so that $c= a_j$ for
some $j\in \{1,\ldots,\ell\}$. But if $n_j(\psi)$ is not prime there exists an integer $1<k<n_j(\psi)$ such that
$k$ divides $n_j(\psi)$. But then, by (3.32), $\g_1$ and $\g[a_j]^k$ generate a proper semisimple subalgebra of
$\g$, contradicting the maximality of $\g_1$.

 \vskip 1.5pc \centerline{\bf 4. Example}\vskip 1pc {\bf 4.1.} In
this section we consider the example of the theory above for the case where, 
for a positive integer $n> 1$, $\g = \hbox{\rm Lie}\,Sl(n,\Bbb C)$. With the usual meaning of matrix units, 
$e_{ij}$, $x\in\g$,  when we can write
$$x =\sum_{i,j =1}^n\,a_{ij}(x)\,\,e_{ij}$$ where $\sum_{i=1}^n \,a_{ii}(x) = 0$. For $k$ a positive
integer,  where $1<k\leq n$,  let $$\delta =\{d_1,\ldots,d_k\}\eqno (4.1)$$ where the $d_p$ are positive integers
such that $$\sum_{p=1}^k\,d_p= n. \eqno (4.2)$$ For $q\in \{1,\ldots,k\}$, put $$f_q = \sum_{p=1}^q\,d_p$$ and
hence
$$1\leq f_1<\cdots < f_k = n. \eqno (4.3)$$ Now for any $i\in \{1,\ldots,n\}$,  let $q(i)\in \{1,\ldots,k\}$ be
the minimum value of $q$ such that $i\leq f_q$. Thus if we put $f_0 = 0$,  and we let $I_q$ be the half-open
interval of integers given by putting $I_q = (f_{q-1},f_q]$,  then one has the disjoint union $$(0,n] =
\sqcup_{q=1}^k I_q.\eqno (4.4)$$ Clearly one has $$\hbox{\rm card}\,I_q = d_q\eqno (4.5)$$ and for 
any $i\in (0,n]$ one has
$$i\in I_{q(i)}. \eqno (4.6)$$ Next put $$\n(\delta) = \{x\in \g\mid a_{ij}(x) =
0,\,\,\hbox{unless}\,\,j>f_{q(i)} 
\}.$$ In addition for $r,s\in \{1,\ldots,q\}$, where $r\neq s$ let $$\g_{r,s}(\delta)=\{x\in \g\mid 
a_{ij}(x) = 0,\,\,\hbox{unless}\,\,i\in I_r,\hbox{and}\,\,j\in I_s\}. $$ One readily notes that $\n(\delta)$ is
a nilpotent Lie algebra and one has the vector space direct sum
$$\n(\delta) = \oplus_{r<s}\,\,\g_{r,s}(\delta). \eqno (4.7)$$ Let $\overline{\n}(\delta)$ be the transpose of
$\n(\delta)$. One then has the direct sum $$\overline{\n}(\delta) = \oplus_{s<r}\,\,\g_{r,s}(\delta).\eqno
(4.8)$$ Also
for
$q\in
\{1,\ldots,k\}$, let
$$\ss_q(\delta) =
\{x\in
\g\mid
 a_{ij}(x)= 0, 
\,\,\hbox{unless}\,\,i,j\in I_q\}.$$ One readily has that $\ss_q(\delta)=0$ if $d_q =1$ and otherwise 
$\ss_q(\delta)$ is a simple
Lie subalgebra of $\g$ where in fact $$\ss_q(\delta)\cong \hbox{\rm Lie}\,Sl(d_q,\Bbb C).\eqno (4.9)$$ Let
$\ss(\delta)$ be the semisimple Lie subalgebra given by putting $$\ss(\delta) = \oplus_{q=1}^k \ss_q(\delta).$$
Now let $\hh$ be the space of all diagonal matrices in $\g$ so that $\hh$ is a Cartan subalgebra of $\g$. Let
$$\tt(\delta) =\{x\in
\hh\mid a_{ii}(x) = a_{jj}(x)\,\,\hbox{if}\,\,q(i) = q(j)\}. $$ Let $$\m(\delta) = \tt(\delta) + \ss(\delta)\eqno
(4.10)$$ and put
$\q(\delta) = \m(\delta) + \n(\delta)$. The following proposition is straightforward and is left as an
exercise.\vs 

{\bf Proposition 4.1.} {\it $\q = \q(\delta)$ is a parabolic subalgebra of $\g$ and, up to
conjugacy, every proper parabolic subalgebra is of this form. Moreover $\q = \m +\n$ is a Levi decomposition of
$\q$ with $\n$ as a nilradical and $\m$ as a Levi factor where $\m = \m(\delta)$ and $\n = \n(\delta)$.
Furthermore (4.10) is the decomposition (1.4) where $\tt = \tt(\delta)$ and $\ss = \ss(\delta)$. Next the set,
$R$, of
$\tt$-roots $\nu$ is parameterized by all pairs $r,s\in \{1,\ldots,k\}$,  where $r\neq s$,  and if the
parameterization is denoted by $\nu(r,s)$, then for any $x\in \tt$ one has $$\nu(r,s)(x) = a_{ii}(x) -
a_{jj}(x)\eqno (4.11)$$ for $i\in I_r$ and $j\in I_s$. In addition  the $\tt$-root space corresponding to
$\nu(r,s)$ is given by $$\g_{\nu(r,s)} = \g_{r,s}(\delta).\eqno (4.12)$$ The irreducible adjoint action of $\m$ 
on $\g_{\nu(r,s)}$ is given as follows: Put 
$\ss_q =
\ss_q(\delta)$. Then $\ss_p$ operates trivially if $p\notin \{r,s\}$. Furthermore
$\g_{\nu(r,s)}$ is one-dimensional if and only if $d_r =d_s =1$. If $d_r =1$ and $d_s >1$ (resp. $d_r >1$ and
$d_s=1$), then
$\g_{\nu(r,s)}$ is $d_s$ (resp. $d_r$)-dimensional and affords a fundamental irreducible $d_s$ (resp.
$d_r$)-dimensional of
$\ss_s$ (resp. $\ss_r$). Moreover if $d_r$ and $d_s$ are both greater than 1,  then
$\hbox{\rm dim}\,\g_{\nu(r,s)}= d_r\,d_s$ and $\g_{\nu(r,s)}$ affords the direct product of a fundamental
irreducible $d_r$-dimensional representation of
$\ss_r$ and a fundamental irreducible $d_s$-dimensional representation of
$\ss_s$. }

\vskip 1.2pc 
\centerline{\bf References}\vskip 1pc
\parindent=42pt

\item {[B-de S]} A. Borel and J. de Siebenthal, Les sous-groupes ferm\'es de rang maximum des Lie clos.
{\it Comment. Math. Helv.}, {\bf 23} (1949), 200--221.
\item {[Hum]} J. Humphreys, {\it Introduction to Lie Algebras and Representation Theory}, Grad. Texts in
Math., {Vol.~9}, Springer-Verlag, New York, 1972. 
\item {[Jac]} N. Jacobson, {\it Lie Algebras}, Wiley(Interscience), New York, {\bf 10},  1962.
\item {[Jos]} A. Joseph, Orbital Varieties of the Minimal orbit, {\it Ann. Ec. Norm. Sup.},{\bf 31}(1998), 
17--45.
\item {[Kos]} B. Kostant, A Characterization of the Classical Groups, {\it Duke Math. J.} {\bf 25}:1(1958),
107--124.
\item {[Wol]} J. Wolf, {\it Spaces of Constant Curvature}, McGraw-Hill, New York,
1967.

\smallskip
\parindent=30pt
\baselineskip=14pt
\vskip 1pc
\vbox to 60pt{\hbox{Bertram Kostant}
      \hbox{Dept. of Math.}
      \hbox{MIT}
      \hbox{Cambridge, MA 02139}}\vskip 2pt

      \noindent E-mail kostant@math.mit.edu

\end